\documentclass{article}
\usepackage{amssymb,amsfonts,amsmath,amsthm,enumerate}
\usepackage[numeric]{amsrefs}
\usepackage{mathtools}

\def \f{\frac}

\def \z{\zeta}

\def \F{{\mathbb F}}
\def \G{\Gamma}

\def \Prim{\mathrm{Prim}}

\def \Re{\operatorname{Re}}

\def \Z{{\mathbb Z}}
\def \l{\left}
\def \r{\right}
\def \hat{\widehat}

\def \z{\zeta}

\newtheorem{theorem}{Theorem}

\numberwithin{equation}{section} \theoremstyle{definition}

\begin{document}

\title{Functional equations for Selberg zeta functions with Tate motives}
\author{Shin-ya Koyama\footnote{Department of Biomedical Engineering, Toyo University,
2100 Kujirai, Kawagoe, Saitama, 350-8585, Japan.} \ \& Nobushige Kurokawa\footnote{Department of Mathematics, Tokyo Institute of Technology, 
Oh-okayama, Meguro-ku, Tokyo, 152-8551, Japan.}}

\maketitle

\begin{abstract}
For a compact Riemann surface $M$ of genus $g\ge 2$, we study the functional equations
of the Selberg zeta functions attached with the Tate motives $f$.
We prove that certain functional equations hold if and only if $f$ has the absolute automorphy.
\end{abstract}



\section*{Introduction}
For a compact Riemann surface $M$ of genus $g\ge2$ the standard Selberg zeta function $Z_M(s)$ is constructed as
$$
Z_M(s)=\prod_{P\in\Prim(M)}\prod_{n=0}^\infty\l(1-N(P)^{-s-n}\r),
$$
where $\Prim(M)$ denotes the set of primitive closed geodesics and the norm $N(P)$ is defined by
$$
N(P)=\exp(\mathrm{length}(P)).
$$
It has the functional equation under $s\rightarrow 1-s$:
$$
Z_M(1-s)=Z_M(s)\exp\l((4-4g)\int_0^{s-\f12}\pi t \tan(\pi t)dt\r).
$$
This functional equation was proved by Selberg \cite{S1, S2} and the following symmetric version was found later:
$$
\hat Z_M(1-s)=\hat Z_M(s),
$$
where
$$
\hat Z_M(s)=Z_M(s)\G_M(s)
$$
with 
$$
\G_M(s)=(\G_2(s)\G_2(s+1))^{2g-2}.
$$
This double gamma function $\G_2(s)$ is the normalized one used in \cite{KK} and we will recall the proof of the
symmetric functional equation for $Z_M(s)$ in the text.

Now, the simple Euler product
$$
\z_M(s)=\prod_{P\in\Prim(M)}(1-N(P)^{-s})^{-1}=\f{Z_M(s+1)}{Z_M(s)}
$$
was also studied later and it is a more natural analog of the usual Euler product for the Riemann zeta function
$$
\z(s)=\prod_{p:\,\text{primes}}(1-p^{-s})^{-1}.
$$
Especially the proof of the prime number theorem
$$
\pi(x)\sim\f x{\log x}\quad(x\to\infty)
$$
applied to $\z_M(s)$ gives the prime geodesic theorem
$$
\pi_M(x)\sim\f x{\log x}\quad(x\to\infty),
$$
where
$$
\pi_M(x)=\#\{P\in\Prim(M)\ |\ N(P)\le x\}.
$$
The functional equation of $\z_M(s)$ has the following form:
$$
\z_M(-s)=\z_M(s)^{-1}(2\sin(\pi s))^{4-4g}.
$$

In this paper we study the functional equations for $\z_{M(f)}(s)$ with Tate motives $f$.
Here we define $\z_{M(f)}(s)$ as
$$
\z_{M(f)}(s)=\prod_k\z_M(s-k)^{a(k)}
$$
for a Laurent polynomial
$$
f(x)=\sum_{k\in\Z}a(k)x^k\in\Z[x,x^{-1}].
$$

It may be suggestive to consider $x=\mathbb T$ the Tate twist.
Of course $\z_{M(1)}(s)=\z_M(s)$ in our notation.

We describe results on $\z_{M(f)}(s)$ only for ``odd'' $f$ here in Introduction.
See the text concerning the ``even'' cases.

\noindent{\bf Theorem}
Let $M$ and $f$ be as above.
For each integer $D$ the following conditions are equivalent.
\begin{enumerate}[\rm (1)]
\item $\z_{M(f)}(D-s)=\z_{M(f)}(s)$.
\item $f(x^{-1})=-x^{-D}f(x)$.
\end{enumerate}

\noindent{\bf Remark}
Condition (2) is called the {\it absolute automorphy} \cite{KT}.
In the paper \cite{KT} the definition of absolute automorphic forms are described in a more general setting
for any function $f$ on positive real numbers, and
the theory of absolute zeta functions $\z_f(s)$ is developed,
which are the autmorphic $L$-functions constructed from the absolute automorphic forms $f$.
It is in the framework of absolute mathematics \cite{DKK, KO}.

For example, let $f(x)=(x-1)^r$ for an odd integer $r\ge1$. Then we see that
$$
f(x^{-1})=-x^{-r}f(x).
$$
Hence Theorem gives the functional equation of $\z_{M(f)}(s)$ as
$$
\z_{M(f)}(r-s)=\z_{M(f)}(s).
$$
A remarkable point is that we need no ``gamma factors'' here.
In the simplest case $r=1$ we get the functional equation for
$$
\z_{M(f)}(s)=\f{\z_{M}(s-1)}{\z_{M}(s)}
$$
as
$$
\z_{M(f)}(1-s)=\z_{M(f)}(s).
$$
We remark that the study of the functional equations for
$$
Z_{M(f)}(s)=\prod_k Z_{M}(s-k)^{a(k)}
$$
is quite similar.

We add a few more comments on $Z_{M(f)}(s)$. Let
$$
f(x)=\sum_k a(k)x^k\in\Z[x,x^{-1}]
$$
satisfying 
$$
f(x^{-1})=Cx^{-D}f(x)
$$
with $C=\pm1$. Then
$$
Z_{M(f)}(s)=\prod_k Z_M(s-k)^{a(k)}
$$
has the functional equation
$$
Z_{M(f)}(D+1-s)=Z_{M(f)}(s)^C S_{M(f)}(s)^C,
$$
where
$$
S_{M(f)}(s)=\prod_k S_{M}(s-k)^{a(k)}
$$
with
$$
S_M(s)=\f{\G_M(s)}{\G_M(1-s)}=(S_2(s)S_2(s+D))^{2-2g}.
$$
Here
$$
S_2(s)=\f{\G_2(2-s)}{\G_2(s)}
$$
is the normalized double sine function of \cite{KK}.
For example $f(x)=x^{-1}-1$ $(C=-1$, $D=-1$) gives the functional equation for
$$
Z_{M(f)}(s)=\f{Z_M(s+1)}{Z_M(s)}=\z_M(s)
$$
as
$$
\z_M(-s)=\z_M(s)^{-1}(2\sin(\pi s))^{4-4g}
$$
where the result
$$
S_{M(f)}(s)=\f{S_M(s+1)}{S_M(s)}=\l(\f{S_2(s+2)}{S_2(s)}\r)^{2-2g}=(2\sin(\pi s))^{4-4g}
$$
is used. Similarly we obtain the functional equation for
$$
Z_{M(f^2)}(s)=\f{Z_M(s+2)Z_M(s)}{Z_M(s+1)}=\f{\z_M(s+1)}{\z_M(s)}=\z_{M(f)}(s)
$$
as
$$
\z_{M(f)}(-1-s)=\z_{M(f)}(s)
$$
that is
$$
Z_{M(f^2)}(-1-s)=Z_{M(f^2)}(s)
$$
with no gamma factors.

\section{Selberg zeta functions}
We describe the needed functional equations for $Z_M(s)$ and $\z_M(s)$ with simple proofs.
Let $\G_r(s)$ be the normalized gamma function of order $r$ defined by
$$
\G_r(s)=\exp\l(\l.\f\partial{\partial w}\z_r(w,s)\r|_{w=0}\r)
$$
with the Hurwitz zeta function of order $r$
$$
\z_r(w,s)=\sum_{n_1,\cdots,n_r\ge0}(n_1+\cdots+n_r+s)^{-w}.
$$
The normalized sine function $S_r(s)$ of order $r$ is constructed as
$$
S_r(s)=\G_r(s)^{-1}\G_r(r-s)^{(-1)^r}:
$$
see \cite{KK} for detailed properties with proofs.

\begin{theorem}
\begin{enumerate}[\rm (1)]
\item Let
$$
\hat Z_M(s)=Z_M(s)\G_M(s)
$$
with
$$
\G_M(s)=(\G_2(s)\G_2(s+1))^{2g-2}.
$$
Then
$$
\hat Z_M(1-s)=\hat Z_M(s).
$$
\item
$$
\z_M(-s)=\z_M(s)^{-1}(2\sin(\pi s))^{4-4g}.
$$
\end{enumerate}
\end{theorem}

{\it Proof.}
(1) From the functional equation for $Z_M(s)$ due to Selberg \cite{S1, S2}
$$
Z_M(1-s)=Z_M(s)\exp\l((4-4g)\int_0^{s-\f12}\pi t \tan(\pi t)dt\r)
$$
we see that it is sufficient to show the identity
$$
\exp\l((4-4g)\int_0^{s-\f12}\pi t \tan(\pi t)dt\r)=\f{\G_M(s)}{\G_M(1-s)}.
$$
We first show that
\begin{equation}
\exp\l((4-4g)\int_0^{s-\f12}\pi t \tan(\pi t)dt\r)
=(S_2(s)S_2(s+1))^{2-2g}.
\end{equation}
Since both sides are equal to 1 at $s=\f12$ (note that
$S_2(\f32)=\G_2(\f32)\G_2(\f12)^{-1}=S_2(\f12)^{-1}$), 
it suffices to show the coincidence of logarithmic derivatives.
The left hand side becomes
$$
(4-4g)\pi\l(s-\f12\r)\tan\l(\pi\l(s-\f12\r)\r)
=(2-2g)\pi(1-2s)\cot(\pi s).
$$
Concerning the right hand side, the differential equation
$$
S_2'(s)=\pi(1-s)\cot(\pi s)S_2(s)
$$
proved in \cite{KK} gives
\begin{align*}
\lefteqn{(2-2g)\l(\f{S_2'(s)}{S_2(s)}+\f{S_2'(s+1)}{S_2(s+1)}\r)}\\
&=(2-2g)\l(\pi(1-s)\cot(\pi s)+\pi(-s)\cot(\pi(s+1))\r)\\
&=(2-2g)\pi(1-2s)\cot(\pi s).
\end{align*}
Thus we obtain (1.1).

Next from (1.1) we get
\begin{align*}
\exp\l((4-4g)\int_0^{s-\f12}\pi t \tan(\pi t)dt\r)
&=(S_2(s)S_2(s+1))^{2-2g}\\
&=\l(\f{\G_2(2-s)}{\G_2(s)}\cdot\f{\G_2(1-s)}{\G_2(s+1)}\r)^{2-2g}\\
&=\f{(\G_2(s)\G_2(s+1))^{2g-2}}{(\G_2(1-s)\G_2(2-s))^{2g-2}}\\
&=\f{\G_M(s)}{\G_M(1-s)}.
\end{align*}
Hence we have the functional equation
$$
Z_M(1-s)=Z_M(s)\f{\G_M(s)}{\G_M(1-s)}
$$
that is 
$$
\hat Z_M(1-s)=\hat Z_M(s)
$$
as desired.

(2) Since
$$
\z_M(s)=\f{Z_M(s+1)}{Z_M(s)}
$$
we have
\begin{align*}
\z_M(-s)\z_M(s)
&=\f{Z_M(1-s)}{Z_M(-s)}\cdot \f{Z_M(s+1)}{Z_M(s)}\\
&=\f{Z_M(1-s)}{Z_M(s)}\cdot \f{Z_M(s+1)}{Z_M(-s)}.
\end{align*}
Hence (1) gives
\begin{align*}
\z_M(-s)\z_M(s)
&=\f{\G_M(s)}{\G_M(1-s)}\cdot \f{\G_M(-s)}{\G_M(s+1)}\\
&=(S_2(s)S_2(s+1))^{2-2g}(S_2(s+1)S_2(s+2))^{2g-2}\\
&=\l(\f{S_2(s+2)}{S_2(s)}\r)^{2g-2}.
\end{align*}
Recall the relations proved in \cite{KK}:
\begin{align*}
S_2(s+2)
&=S_2(s+1)S_1(s+1)^{-1}\\
&=S_2(s+1)(-2\sin(\pi s))^{-1}
\end{align*}
and
\begin{align*}
S_2(s)
&=S_2(s+1)S_1(s)\\
&=S_2(s+1)(2\sin(\pi s)).
\end{align*}
Thus we get the functional equation for $\z_M(s)$:
$$
\z_M(-s)\z_M(s)
=(2\sin(\pi s))^{4-4g}
$$
that is
$$
\z_M(-s)=\z_M(s)^{-1}(2\sin(\pi s))^{4-4g}.
$$
\hfill\qed

\section{Functional equation for $\z_{M(f)}(s)$}
Let
$$
\z_{M(f)}(s)=\prod_k\z_M(s-k)^{a(k)}
$$
for
$$
f(x)=\sum_{k}a(k)x^k\in\Z[x,x^{-1}].
$$
We prove the following theorem.

\begin{theorem}
For each integer $D$ the following conditions are equivalent:
\begin{enumerate}[\rm(1)]
\item $\z_{M(f)}(D-s)=\z_{M(f)}(s)$.
\item $f(x^{-1})=-x^{-D}f(x)$.
\item $a(D-k)=-a(k)$ for all $k$.
\end{enumerate}
\end{theorem}

{\it Proof.}
We first show the equivalence $(2)\Longleftrightarrow(3)$.
Let
$$
f(x)=\sum_{k}a(k)x^k.
$$
Then
\begin{align*}
x^D f(x^{-1})
&=\sum_k a(k)x^{D-k}\\
&=\sum_k a(D-k)x^k,
\end{align*}
where we needed the exchange $k\longleftrightarrow D-k$. Hence
$$
x^Df(x^{-1})=-f(x)
$$
is equivalent to
$$
a(D-k)=-a(k)\qquad \text{for all }k.
$$
Next we show the equivalence $(1)\Longleftrightarrow(2)$. Since
\begin{align*}
\z_{M(f)}(D-s)
&=\prod_k\z_M((D-s)-k)^{a(k)}\\
&=\prod_k\z_M((D-k)-s)^{a(k)}\\
&=\prod_k\z_M(k-s)^{a(D-k)},
\end{align*}
the functional equation for $\z_M(s)$ gives
\begin{align*}
\z_{M(f)}(D-s)
&=\prod_k (\z_M(s-k)^{-1}(2\sin(\pi s))^{4-4g})^{a(D-k)}\\
&=\l(\prod_k\z_M(s-k)^{-a(D-k)}\r)(2\sin(\pi s))^{(4-4g)f(1)},
\end{align*}
where we used
$$
f(1)=\sum_k a(k)=\sum_k a(D-k).
$$
Hence we have the following expression
\begin{equation}
\f{\z_{M(f)}(D-s)}{\z_{M(f)}(s)}
=\l(\prod_k\z_M(s-k)^{-a(D-k)-a(k)}\r)(2\sin(\pi s))^{(4-4g)f(1)}.
\end{equation}
From this expression the equivalence $(1)\Longleftrightarrow(3)$ is shown as follows.
First the condition (3) (or equivalently (2)) implies $f(1)=0$ and that
$a(D-k)+a(k)=0$ for all $k$. Hence (2.1) gives
$$
\f{\z_{M(f)}(D-s)}{\z_{M(f)}(s)}=1,
$$
which is (1).

Now assume (1). Then from (2.1) we have the identity
\begin{equation}
\prod_k\z_M(s-k)^{a(D-k)+a(k)}=(2\sin(\pi s))^{(4-4g)f(1)}.
\end{equation}
Since $\z_M(s)$ is non-zero holomorphic in $\Re(s)>1$, the left hand side of (2.2) is
non-zero holomorphic for sufficiently large $\Re(s)$. Hence looking at the left hand side at
sufficiently large $s\in\Z$ we see $f(1)=0$. Then (2.2) gives
\begin{equation}
\prod_k\z_M(s-k)^{a(D-k)+a(k)}=1.
\end{equation}
We remark that (2.3) is actually written as
\begin{equation}
\prod_{k\le K}\z_M(s-k)^{a(D-k)+a(k)}=1
\end{equation}
for some $K\in\Z$, since $f(x)\in\Z[x,x^{-1}]$. Hence we have the identity
\begin{equation}
\z_M(s-K)^{a(D-K)+a(K)}=\prod_{k< K}\z_M(s-k)^{-a(D-k)-a(k)}.
\end{equation}
Look at (2.5) at $s=K+1$, then the right hand side is 
$$
\prod_{k< K}\z_M(1+(K-k))^{-a(D-k)-a(k)},
$$
which is a finite non-zero value. Hence looking at the left hand side of (2.5) we see that
$$
a(D-K)+a(K)=0
$$
since $\z_M(s-K)$ has a simple pole at $s=K+1$. Thus (2.5) becomes
\begin{equation}
\prod_{k\le K-1}\z_M(s-k)^{a(D-k)+a(k)}=1
\end{equation}
Inductively we see (3).
\hfill\qed

Theorem 2 treated ``odd'' $f$.
The next theorem deals with the other case for ``even'' $f$.

\begin{theorem}
For each integer $D$ the following conditions are equivalent:
\begin{enumerate}[\rm(1)]
\item $\z_{M(f)}(D-s)=\z_{M(f)}(s)^{-1}(2\sin(\pi s))^{(4-4g)f(1)}$.
\item $f(x^{-1})=x^{-D}f(x)$.
\item $a(D-k)=a(k)$ for all $k$.
\end{enumerate}
\end{theorem}

{\it Proof.}
The equivalence $(2)\Longleftrightarrow(3)$ is shown exactly as in the proof of Theorem 2.
Now we show $(1)\Longleftrightarrow(3)$. Notice that
\begin{align*}
\z_{M(f)}(D-s)
&=\prod_k\z_M((D-s)-k)^{a(k)}\\
&=\prod_k\z_M((D-k)-s)^{a(k)}\\
&=\prod_k\z_M(k-s)^{a(D-k)}\\
&=\prod_k\l(\z_{M}(s-k)^{-1}(2\sin(\pi s))^{4-4g}\r)^{a(D-k)}\\
&=\l(\prod_k \z_{M}(s-k)^{-a(D-k)}\r)(2\sin(\pi s))^{(4-4g)f(1)},
\end{align*}
where we used that
$$
\sum_k a(D-k)=f(1).
$$

\underline{\it Proof of $(3)\Longrightarrow(1)$.}
From (3) we have
\begin{align*}
\z_{M(f)}(D-s)
&=\l(\prod_k \z_M(s-k)^{-a(k)}\r)(2\sin(\pi s))^{(4-4g)f(1)}\\
&=\z_{M(f)}(s)^{-1}(2\sin(\pi s))^{(4-4g)f(1)},
\end{align*}
which is (1).

\underline{\it Proof of $(1)\Longrightarrow(3)$.}
Since
$$
\z_{M(f)}(D-s)=\l(\prod_k \z_M(s-k)^{-a(D-k)}\r)(2\sin(\pi s))^{(4-4g)f(1)}
$$
as above, we have
$$
\f{\z_{M(f)}(D-s)}{\z_{M(f)}(s)^{-1}(2\sin(\pi s))^{(4-4g)f(1)}}
=\prod_k \z_M(s-k)^{a(k)-a(D-k)}.
$$
Hence from the assumption (1) we get
$$
\prod_k \z_M(s-k)^{a(k)-a(D-k)}=1,
$$
which can be written as
$$
\prod_{k\le K} \z_M(s-k)^{a(k)-a(D-k)}=1
$$
that is
$$
\z_M(s-K)^{a(K)-a(D-K)}=\prod_{k< K} \z_M(s-k)^{a(D-k)-a(k)}.
$$
Then we obtain $a(D-K)=a(K)$ and inductively $a(D-k)=a(k)$ for all $k$ exactly as in the proof of Theorem 2.
\hfill\qed

{\bf Example.}
Let $f(x)=(x-1)^r$ for an even integer $r\ge0$. Then we see that
$$
f(x^{-1})=x^{-r}f(x).
$$
Hence we obtain the functional equation
$$
\z_{M(f)}(r-s)=\z_{M(f)}(s)^{-1}\times
\begin{cases}
(2\sin(\pi s))^{4-4g}&(r=0),\\
1&(r\ge2,\text{ even}).
\end{cases}
$$
Of course the $r=0$ case gives the functional equation of $\z_M(s)$.

{\bf Remark.}
Let $f(x)=(x-1)^r$ for an integer $r\ge0$. Then $\z_M(s)$ is written explicitly as
$$
\z_{M(f)}(s)=\prod_{k=0}^r \z_M(s-k)^{(-1)^{r-k}\binom rk}.
$$
In this case another suggestive notation would be
$$
\z_{M(f)}(s)=\z_{M\otimes\mathbb G_m^r}(s)
$$
since $(x-1)^r$ is the counting function of $\mathbb G_m^r$;
see \cite{CC, KO, DKK}.

\end{document}